\documentclass[12pt]{amsart}
\usepackage{amsmath,amssymb}
\usepackage[mathscr]{eucal}
\begin{document}

\newcommand{\bbR}{{\mathbb R}}
\newcommand{\bbN}{{\mathbb N}}
\newcommand{\bbI}{{\mathbb I}}
\newcommand{\bbQ}{{\mathbb Q}}
\newcommand{\bbP}{{\mathbb P}}
\def\A{\mathscr{A}}
\def\B{\mathscr{B}}
\def\C{\mathscr{C}}
\def\D{\mathscr{D}}
\def\E{\mathscr{E}}
\def\F{\mathscr{F}}
\def\H{\mathscr{H}}
\def\K{\mathscr{K}}
\def\M{\mathscr{M}}
\def\N{\mathscr{N}}
\def\G{\mathscr{G}}
\def\P{\mathscr{P}}
\def\U{\mathscr{U}}
\def\V{\mathscr{V}}
\def\W{\mathscr{W}}
\renewcommand{\ss}{\subseteq}
\newcommand{\dsize}{\displaystyle}
\def\dom{\mathop{\mathrm{dom}}}
\newtheorem{thm}{Theorem}[section]
\newtheorem{cor}[thm]{Corollary}
\newtheorem{lem}[thm]{Lemma}
\newtheorem{prop}[thm]{Proposition}
\newtheorem{remark}[thm]{Remark}
\newtheorem{ex}[thm]{Example}
\font\sans=cmss10
\def\axiom#1{{\sans (#1)}}
\def\wk{\text{weak}}
\def\domc#1{{\dom c_{#1}}}
\def\sm{\smallsetminus}

\title{On closed mappings of $\sigma$-compact spaces and  dimension}
\author{El\.{z}bieta Pol and Roman Pol}
%\thanks{The research of the second
%author was partially supported by KBN grant 2PO3A01724}
\address{University of Warsaw \and  University of Warsaw}
\email{E.Pol@mimuw.edu.pl \and  R.Pol@mimuw.edu.pl}
\date{\today}

\begin{abstract}
A remainder of the Hilbert space $l_{2}$ is a space homeomorphic to $Z \setminus l_{2}$, where $Z$ is a metrizable compact extension of $l_{2}$, with $l_{2}$ dense in $Z$. We prove that for any remainder $K$ of $l_{2}$, every non-one-point closed image of $K$  either contains a compact set with no transfinite dimension  or contains compact sets of arbitrarily high inductive transfinite dimension ind. We shall also construct for each natural  $n$ a $\sigma$-compact metrizable $n$-dimensional space whose image under any non-constant closed map has dimension at least $n$, and analogous examples for the transfinite inductive dimension ind (this provides a rather strong negative solution of a problem in \cite{EnPo}).
\end{abstract}

\maketitle

\section{Introduction}

We shall consider only   metrizable separable spaces. %Our terminology follows \cite{Ke} and \cite{E1).
A remainder of a space $E$ is a space homeomorphic to $\tilde{E} \setminus E$, where $\tilde{E}$ is a compact extension of $E$ with $E$ dense in $\tilde{E}$, cf. \cite{AN}, \cite{E1}.
Notice that if $E$ is topologically complete, its remainders are countable unions of compact sets.

A standard remainder of the Hilbert space $l_{2}$ is the subspace $K_{\omega}$ of the Hilbert cube $[0,1] ^{\bbN}$ consisting of points with all but finitely many coordinates zero, cf. \cite{BP}. However, remainders of $l_{2}$ form a much wider class, cf. \cite{B}, \cite{C},  \cite{HeW}, \cite{vM},  including spaces which are countable unions of finite-dimensional compacta
containing no arcs, cf. \cite{vM}.

Let us recall that  the transfinite dimension ind is the extension of the Menger-Urysohn dimension by transfinite induction, cf. \cite{HuW} (some additional information is in section 2).

The space $K_{\omega}$ has no transfinite dimension and in fact, by \cite{EP}, no image of $K_{\omega}$ under a closed map has transfinite dimension, unless it reduces to a singleton.

We shall strengthen this result to the following effect.

\begin{thm} Let $K$ be a remainder of the Hilbert space $l_{2}$. Then every  non-one-point image of $K$ under a closed map  either contains a compact set with no transfinite dimension  or
contains compact sets of arbitrarily high transfinite dimension.
\end{thm}

The proof is based  on \cite{EP} and our recent paper \cite{PP} which appeals to the classical theorem of Hurewicz on non-analyticity of the collection of compact sets of the rationals.

Similar reasonings, although in a different setting, yield the following result, which provides a solution (in a rather strong form) to Problem 7.19 in \cite{EnPo}.

\begin{thm} For each countable ordinal $\alpha \geq 2$,  there exists a $\sigma$-compact space $Z$ which has the transfinite dimension ind (if $\alpha$ is finite, so is ind$Z$), such that every image of $Z$ under a non-constant closed mapping contains a compact set whose transfinite dimension ind is
 either $\geq \alpha$ or is not defined.
\end{thm}

The $\sigma$-compact spaces $Z$ in Theorem 1.2 can not be completely metrizable. Indeed, if a non-one-point space contains a compact set with nonempty interior, then it admits a non-constant closed map into the real line.

However, if we restrict the class of mappings to perfect maps, the situation is not  clear. We address this briefly in Comments 5.4 and 5.4.

\section{Some background}

\subsection{The transfinite dimension ind and countable-dimensional spaces.}
The transfinite dimension ind is the extension by transfinite induction of the Menger-Urysohn small inductive dimension, see \cite{HuW}, \cite{Na}, \cite{E2}.
All spaces for which the transfinite inductive dimension ind is defined are countable-dimensional, i.e., they are  countable unions of zero-dimensional sets.
However, the space $K_{\omega}$ is countable-dimensional, but has no transfinite dimension
 ($K_{\omega}$ contains topologically all Smirnov's spaces $\;S_{\alpha}\,$, $\;\alpha < \omega _{1}\,$,  and the set  $\{ {\rm ind} S_{\alpha}: \;\alpha < \omega _{1} \}$ is unbounded in $\omega _{1}$, cf. \cite{E2}, 5.3.11,   7.1.33 and 7.2.12).

 Let us notice also that a space $E$ has the transfinite dimension ind  if and only if $E$ has a countable-dimensional topologically complete (equivalently - compact) extension, cf. \cite{E2}, 7.2.19 and 7.2.21.

We shall need an observation (easily checked by  induction), that if $F$ is a closed nonempty subspace of a space $X$ and ind$(X \setminus F) = 0$ then ind$X = {\rm ind}F$.

We shall also use the fact that perfect mappings with zero-dimensional fibers do not lower the transfinite dimension ind.

\subsection{The Effros Borel spaces.}
%Our terminology related to the descriptive set theory follows \cite{D} and \cite{K}.

Given a space $E$, we denote by $F(E)$ the space of closed subsets of $E$ and   the Effros Borel structure in $F(E)$
 is the $\sigma$-algebra in $F(E)$ generated by the sets $\{ A \in F(E): A \cap U \neq \emptyset \}$, where $U$ is open in $E$. %, cf. \cite{Ke}, 12.C.
 The Effros Borel space $F(E)$ is the space $F(E)$  equipped with the Effros Borel structure  and if $E$ is either $\sigma$-compact or completely metrizable the
 Effros Borel space $F(E)$ is standard, cf. \cite{K}, 12.C.
 If $E$ is  compact the Effros Borel structure coincides with the $\sigma$-algebra of Borel sets with respect to the Vietoris topology in the hyperspace $F(E)$, cf. \cite{K}, 12.7.

  %If $X$ is a compact metrizable extension of an analytic set $E \subset X$, the map $A \to \overline{A}$ (the closure is taken in $X$) from $F(E)$ to $F(X)$ is a Borel isomorphism, with respect to %the Effros Borel structures, onto the analytic subspace $\{ \overline{A} : \, A \in F(E) \}$ of the hyperspace $F(X)$ and hence Souslin sets in $F(E)$ are mapped onto analytic sets in $F(X)$, cf. %\cite{D}, Section 2. In particular, if $E \subset G \subset X$ and $G$ is analytic, the collection of closures of elements of $F(E)$ in $G$ is a Souslin set in $F(G)$.

\subsection{Non-trivial closed maps of remainders of $l_{2}$ are perfect.}

It was noticed in \cite{EP}, Lemma 3, that if $f : K_{\omega} \to L$ is a closed map onto a non-one-point space, then all fibers of $f$ are compact, i.e., $f$ is perfect. The same is true for any closed map $f : K \to L$ on a remainder $K$ of $l_{2}$, cf. \cite{Mi}.

Let us outline a justification of this fact (being a minor modification of a reasoning from \cite{EP}).

Let $Z$ be a compactum containing $K$, where $Z \setminus K$ is a dense subspace of $Z$ homeomorphic to $l_{2}$. By Va\v{\i}n\v{s}te\v{\i}n's Lemma (cf. \cite{E1}, 4.4.16), for each $x \in L$ the set ${\rm Fr} f^{-1}(x)$ is compact. To prove that $f$ is perfect it suffices to show that ${\rm Int} f^{-1}(x)= \emptyset$ for every $x \in L$. Suppose that Int$f^{-1}(x_{0})\neq \emptyset$ for some $x_{0}\in L$, and let $y_{0} \in {\rm Int}f^{-1}(x_{0})$. Since $L \neq \{ x_{0}\}$, $P= {\rm Fr}f^{-1}(x_{0})$ is a partition in $K$ between $y_{0}$ and some other  point in $K$. This is, however, impossible, since no compact subset of $K$ separates $K$.
To see this, suppose that $P$ is a compact partition in $K$ between points $x$ and $y$. Then there exists a partition $P^{\star}$ in $Z$ between $x$ and $y$ such that $K \cap P^{\star} \subset P$ (cf. \cite{E2}, Lemma 1.2.9). Since $Z \setminus K$ is dense in $Z$, $P^{\star} \cap (Z \setminus K) = P^{\star} \setminus P$ is a $\sigma$-compact partition in $Z \setminus K$ between some two points, which contradicts the fact that no $\sigma$-compact subset of $l_{2}$ separates $l_{2}$ (cf. \cite{BP}, Ch.V, Corollary 6.2 and Theorem 6.3).

\subsection{Hereditarily disconnected $G_{\delta}$-sets in compacta.}
Let $K$ be a compact space with ind$K = \alpha +2$, $\alpha < \omega _{1}$. Then there exists a hereditarily disconnected $G_{\delta}$-set $W$ in $K$ with ind$W \geq \alpha + 1$.

For finite $\alpha$, this can be derived from Theorems 3.9.3 and 3.11.8 in \cite{vM2}.

For arbitrary countable ordinal $\alpha$, one can use the following reasoning from \cite{P2}. Let $d$ be a metric on $K$. There is a point $p$ in $K$ and $0 < r < 1$ such that for each $t \in (0,r)$, the sphere $S_{t} = \{ x : \; d(x,p) = t \}$ has dimension ind$S_{t} \geq \alpha + 1$.

Let $K$ be embedded in the Hilbert cube $I^{\bbN}$, $I = [0,1]$. Then the graph $L = \{ (x,d(x,p)): x \in K \} \subset I^{\bbN} \times I$ is a topological copy  of $K$ and each section $L \cap (I^{\bbN} \times \{ t \})$ is a copy of $S_{t}$. The reasoning in section 6.2 of \cite{P2} provides a $G_{\delta}$-set $W$ in $L$ which projects onto the second coordinate onto the irrationals in $(0,r)$ in a one-to-one way, and ind$W \geq \alpha +1$.

In particular, using Smirnov's compact spaces (cf. \cite{E2}), one can define completely metrizable spaces $D_{\alpha}$ with ind$D_{\alpha} = \alpha$ containing no non-trivial continuum (cf. \cite{P2}, 6.2).

\subsection{An application of Hurewicz's theorem.}

We shall use the following result from \cite{PP}, %Proposition 3.1
based on a classical theorem of Hurewicz, cf. \cite{K}, Exercises 27.8, 27.9.
\begin{prop}
Let $p : E \to F$ be a continuous map of a complete space $E$ onto a non-$\sigma$-compact  space $F$. Then  there is a collection $\D$ of closed relatively discrete sets in $E$
such that for any analytic collection $\A$ in the Effros Borel space $F(E)$ containing $\D$,
 there are $A \in \A$ and $t \in F$ with $p^{-1}(t) \subset A$.
\end{prop}

\section{Proof of Theorem 1.1}

The Hilbert space $l_{2}$ is homeomorphic to the countable product of the  real line $\bbR ^{\bbN}$, cf. \cite{BP}.  Let $K = Z \setminus \bbR ^{\bbN}$,
where $Z$ is a compact extension of $\bbR ^{\bbN}$. By 2.3 it is enough to assume that
\smallskip\begin{enumerate}
\item[(1)] $\; f: K \to L\;$ is a perfect surjection.
\end{enumerate}
We shall also assume that $L \subset [0,1]^{\bbN}$ and using a theorem of Va\v{\i}n\v{s}te\v{\i}n \cite{E1}, Problem 4.5.13(d), one can find a $G_{\delta}$-set $G$ in $Z$ containing $K$ and an extension $\tilde{f}$
of $f$ over $G$ with values in $[0,1]^{\bbN}$ such that
\smallskip\begin{enumerate}
\item[(2)] $\; \tilde{f} : G \to \tilde{f}(G) \;$ is perfect.
\end{enumerate}
Let $Z \setminus G = \bigcup _{n} Z_{n}$, where $Z_{n} \subset \bbR ^{\bbN}$ are compact. The projection of $Z_{n}$ onto the $n'th$ coordinate is contained in an interval $[a_{n},b_{n}]$ and let $c_{n} > b_{n}$.  Then, cf. \cite{BP}, Ch.VI, Example 8.1,
\smallskip\begin{enumerate}
\item[(3)] $\; E = \prod _{n=1}^{\infty} [c_{n}, +\infty )$ is homeomorphic to $\bbR ^{\bbN}$,
\end{enumerate}
and
\smallskip\begin{enumerate}
\item[(4)] $\; \overline{E} \subset G$,
\end{enumerate}
the closure being considered in $Z$.

Let, cf. (4),
\smallskip\begin{enumerate}
\item[(5)] $\; M = \overline{E} \setminus E = K \cap \overline{E}$.
\end{enumerate}

Since $\overline{E}$ is compact and no compact subset of $E$ has nonempty interior in $E$, cf. (3), $M$ is dense in $\overline{E}$.
Since $f \mid M$ is perfect (see (1) and (5)) and $\tilde{f} \mid \overline{E} : \overline{E} \to \tilde{f}(\overline{E} )$ is an extension of $f \mid M$, we infer that (cf. \cite{E1}, Lemma 3.7.4)
\smallskip\begin{enumerate}
\item[(6)] $\; \tilde{f}(M) \cap \tilde{f}(E) = \emptyset$.
\end{enumerate}

Since   $\bbR ^{\bbN}$  is homeomorphic to its square, applying Proposition 2.1 to the projection $p : \bbR ^{\bbN} \times \bbR ^{\bbN} \to \bbR ^{\bbN}$,
we get, by (3),
  a collection $\D$ of closed relatively discrete sets in $E$ such that for any analytic collection $\A$ in the Effros Borel space $F(E)$ containing $\D$, there is an element of $\A$ containing a closed copy of $\bbR ^{\bbN}$.

Let, cf. (3) and (4),
\smallskip\begin{enumerate}
\item[(7)] $\; T= \tilde{f}(E), \;\; S = \overline{T}$.
\end{enumerate}
Then, cf. (1), (2), (5), (6),
\smallskip\begin{enumerate}
\item[(8)] $\; S= \tilde{f}(\overline{E})= T \cup f(M) \subset T \cup L$.
\end{enumerate}
Aiming at a contradiction, assume that all compact subsets of $f(M)$ have transfinite dimension and there exists $\alpha < \omega _{1}$ which bounds the  transfinite dimension ind of all compact sets in $f(M)$.

Since the restriction of $\tilde{f}$ to $E$ is a closed map onto $T$, cf. (6), (7), the collection
\smallskip\begin{enumerate}
\item[(9)] $\;  \tilde{f}(\D) =  \{ \tilde{f}(D): \; D \in \D \} \subset F(T) $
\end{enumerate}
consists of closed, topologically discrete subspaces of $T$. Therefore, for any $B \in \tilde{f}(\D)$, the closure  $\overline{B}$ of $B$ in $S$ intersects $f(M)$ in a compact set, and in effect
\smallskip\begin{enumerate}
\item[(10)] $\;  {\rm ind} \overline{B} \leq \alpha \,$ for any $\,B \in \tilde{f}(\D)  $.
\end{enumerate}
The set of compacta in S with transfinite dimension ind not greater than $\alpha$ is analytic in the hyperspace $F(S)$, cf. \cite{P1}, and the map $\, B \to \overline{B}\,$ from the Effros Borel space $F(T)$ to $F(S)$ being Borel, we conclude that
\smallskip\begin{enumerate}
\item[(11)] $\;  \E = \{ B \in F(T): \; {\rm ind} \overline{B} \leq \alpha \}\;$ is analytic in $F(T)$.
\end{enumerate}
By (9), $\tilde{f}(\D) \subset \E \,$ and, since the map $\, A \to \tilde{f} (A)\,$ from $F(E)$ to $F(T)$ is Borel, the set
\smallskip\begin{enumerate}
\item[(12)] $\;  \A = \{ A \in F(E): \; \hbox{the closure of } \tilde{f}(A) \hbox{ in } S \hbox{ has ind } \leq \alpha \}\;$ is analytic in $F(E)$,
\end{enumerate}
and
\smallskip\begin{enumerate}
\item[(13)] $\;  \D \subset \A $.
\end{enumerate}
By the choice of $\D$, there is $A \in \A$ containing a closed copy $C$ of $\bbR ^{\bbN}$. In effect, we get a perfect map $\; \tilde{f} \mid C : C \to \tilde{f}(C)\;$ of  a copy of $\bbR ^{\bbN}$
onto a space with transfinite dimension ind $ \leq \alpha $, which contradicts \cite{EP}, Remark 2.

This completes the proof of Theorem 1.1.

\section{Proof of Theorem 1.2}

We shall first describe a completely metrizable space $T$ with ind$T = \alpha$, which contains no non-trivial continuum, can not be separated by any closed $\sigma$-compact set and contains a closed set $H$ admitting a continuous map  $h : H \to P$ onto the irrationals such that
\smallskip\begin{enumerate}
\item[(1)] $\;$ ind$h^{-1}(t) = \alpha$ for all $t \in P$.
\end{enumerate}

We begin with a 1-dimensional completely metrizable, connected space $S$ containing no non-trivial continuum \cite{KS}.

 In particular, compact sets in $S$ are boundary, and we can find a discrete collection $P_{1}, P_{2}, \ldots$ of closed copies of the irrationals in $S$.

Let $P$ be the set of irrationals in $I = [0,1]$, $I \setminus P = \{ q_{1}, q_{2}, \ldots \}$ and $S' = (P \times S) \cup \bigcup _{i=1}^{\infty}( \{ q_{i} \} \times P_{i})$.

One easily checks that $S'$  is  1-dimensional, completely metrizable and contains no non-trivial continuum. What we gain by passing from $S$ to $S'$ is that no closed $\sigma$-compact set separates $S'$.

Indeed, aiming at a contradiction, assume that $L \subset S'$ is $\sigma$-compact and $S' \setminus L = U_{0} \cup U_{1}$, where $U_{i}$ are nonempty open sets and $U_{0} \cap U_{1} = \emptyset$.

Let $p : S' \to I$ be the projection. The sets $M_{n} = L \cap (I \times P_{n})$, $n = 1,2,\ldots$, and $M_{0} = L \setminus \bigcup _{n=1}^{\infty} M_{n}$ are $F_{\sigma}$-sets in $L$ and hence $p(M_{i})$ are $\sigma$-compact boundary sets in $I$, for $i = 0, 1,2,\ldots $ Therefore, $A = P \setminus \bigcup _{i=0}^{\infty} p(M_{i})$ is dense in $I$ and $(A \times S) \cap L = \emptyset$. The sets
$A_{i} = \{ a \in A: \, \{a\} \times S \subset U_{i} \}$, $i = 0,1$, cover $A$, $S$ being connected, hence $\overline{A_{0}} \cup \overline{A_{1}} = I$ and we can pick $b \in \overline{A_{0}} \cap \overline{A_{1}}$.  Since $p^{-1}(b)$ is not $\sigma$-compact, there is $c \in p^{-1}(b) \setminus L$. But then $c \in \overline{U_{0}} \cap \overline{U_{1}}$ and $c \not\in L$, which is impossible.

To get the space $T$, we take a completely metrizable space $D_{\alpha}$ with ind$D_{\alpha} = \alpha$, containing no non-trivial continuum, cf. 2.4, and we replace $\{q_{1}\} \times P_{1}$ by a closed copy of the product $P \times D_{\alpha}$. More precisely, we fix a perfect surjection $u : \{q_{1}\} \times P_{1} \to P \times D_{\alpha}$ and we let $T = S' \cup _{u} (P \times D_{\alpha})$ be the result of attaching $P \times D_{\alpha}$ to $S'$ through the map $u$.

Then ind$T =$ ind$(P \times D_{\alpha}) = \alpha$, and since $S'$ maps perfectly onto $T$, no closed $\sigma$-compact set separates $T$. Also, $P \times D_{\alpha}$ embeds in $T$ as a closed subspace $H$. Let $h : P \times D_{\alpha} \to P$ be the projection.

Having defined $T$ and $H$, let us consider a countable-dimensional compact extension $T^{\star}$ of $T$ such that $T$ is dense in $T^{\star}$ (if $\alpha < \omega$, we can have ind$T^{\star} = $ ind$T$), and let
\smallskip\begin{enumerate}
\item[(2)] $\;$ $Z = T^{\star} \setminus T$, $\;K= \overline{H}\;$ and $\;C = K \setminus H$,
\end{enumerate}
where the closure is taken in $T^{\star}$.

Then the arguments from 2.3 show that all closed non-constant maps on $Z$ are in fact perfect.

We shall check that the space $Z$ has the required properties. It is enough to show that if
\smallskip\begin{enumerate}
\item[(3)] $\;$  $f:C \to L\;$ is a perfect surjection,
\end{enumerate}
then for some compact set $B$ in $L$, either $B$ has no transfinite dimension ind or else ind$B \geq \alpha$.

Aiming at a contradiction, assume that
\smallskip\begin{enumerate}
\item[(4)] $\;$  ind$B < \alpha \;$ for any compact set $B \subset L$.
\end{enumerate}

As in section 3, we consider $L \subset I^{\bbN}$, $I = [0,1]$, and we use Va\v{\i}n\v{s}te\v{\i}n's theorem to get a $G_{\delta}$-set $G$ in $K$ containing $C$ and an extension $\tilde{f}$ of $f$ with values in $I^{\bbN}$ such that
\smallskip\begin{enumerate}
\item[(5)] $\;$  $\tilde{f} : G \to \tilde{f}(G)\;$ is  perfect.
\end{enumerate}
The set $K \setminus G \subset H$ is  $\sigma$-compact and $h(K \setminus G) \subset P$. Therefore, there exists a closed in $P$ copy $J$ of the irrationals, disjoint from the $\sigma$-compact set $h(K \setminus G)$. Then $h^{-1}(J)$ is closed in $H$ and hence its closure $\overline{h^{-1}(J)}$ in $K$ is contained in $G$.
Let
\smallskip\begin{enumerate}
\item[(6)] $\;$ $M = \overline{h^{-1}(J)} \cap C$, $\;E = \overline{M} \cap h^{-1}(J)$.
\end{enumerate}
Since $h^{-1}(J) \setminus \overline{M}$ is a $\sigma$-compact subset of a hereditarily disconnected set $H$, it is zero-dimensional. Applying to the fibers $h^{-1}(t)$ an observation from 2.1, and using (1), we infer that, cf. (6),
\smallskip\begin{enumerate}
\item[(7)] $\;$ ind$(h^{-1}(t) \cap E) = \alpha\;$ for $\;t \in J$.
\end{enumerate}
Since the restriction $\tilde{f} \mid M = f \mid M$ is perfect, and $M$ is dense in $E$, we have also
\smallskip\begin{enumerate}
\item[(8)] $\;$ $\tilde{f}(M) \cap \tilde{f}(E)= \emptyset$.
\end{enumerate}

We shall apply Proposition 2.1 to the map $p = h \mid E$ and let $\D$ be a collection of closed, relatively discrete sets in $E$, provided by this proposition.

We can now repeat reasonings from section 3, involving collections $\D$ and $\A$ in the Effros Borel space $F(E)$ almost verbatim, changing in formulas (10), (11) and (12) in section 3 the inequality $\leq \alpha$ to the strict inequality $< \alpha$, up to the point, where one considers an element $A \in \A$ such that for some $t \in J$, $\,p^{-1}(t)  \subset A$. Now, $\tilde{f}$ restricted to $F = p^{-1}(t)$ is perfect, cf. (8), and since $F$ is hereditarily disconnected, the fibers of $\tilde{f} \mid F$ are zero-dimensional. As we noted in sec. 2.1, this implies that ind$\tilde{f}(A) \geq$ ind$\tilde{f}(F) \geq$ ind$F = \alpha$, cf. (7), providing a contradiction with (4), which completes the proof.

\section{Comments}

\subsection{Completions of perfect images of remainders of $l_{2}$.}

Let $K$ be a remainder of $l_{2}$ and let $f : K \to L$ be a perfect surjection. Then every completion $L ^{\star }$ of $L$ contains a strongly infinite-dimensional compactum (cf. \cite{E1}, Ch.6).

Indeed,  the reasoning in section 3 yields a closed copy $\;E\;$ of $\;l_{2}\;$ and a perfect map $\tilde{f} : \overline{E} \to L ^{\star }\;$ ($\overline{E}$ is the closure of $E$ in a compactification $l_{2} \cup K$ with the remainder $K$) such that $\,\tilde{f} (E) \cap \tilde{f} (\overline{E} \setminus E) = \emptyset$. Then $\,\tilde{f} \mid E : E \to L ^{\star } \setminus L$ is a perfect map and \cite{EP}, Remark 2, provides a strongly infinite dimensional compactum in $L ^{\star }$.

\subsection{The transfinite dimension dim.}

 Using similar reasonings as in Section 3 one can show that for any perfect image $L$ of a remainder $K$ of $l_{2}$, either $L$ contains a strongly infinite-dimensional compactum or the transfinite dimension dim (cf. \cite{E2}, 7.3.19) of compacta in $L$ is unbounded.

\subsection{A question.}

There is a countable-dimensional space $E$ which is an absolute $G_{\delta \sigma}$-set, all compact subsets of $E$ are at most zero-dimensional and $E$ has no transfinite dimension ind.

Indeed, there exists a $G_{\delta \sigma}$-set $E$ in $K_{\omega}$ such that neither $E$ nor $K_{\omega} \setminus E$ contain a non-trivial continuum, cf. \cite{PP2}, Remark 4.2.
Then, for any $G_{\delta}$-set $G$ in $K_{\omega}$ containing $E$, $K_{\omega} \setminus G$ is zero-dimensional, and since $K_{\omega}$ has no transfinite dimension, this is also true for $G$.
This shows that $E$ has no transfinite dimension, cf. 2.1.

However, we do not know if there exist countable-dimensional $\sigma$-compact sets $H$ such that $H$ has no transfinite dimension ind but the transfinite dimensions of all compact sets in $H$ are bounded.

\subsection{Perfect maps.}

A part of the reasoning in section 4 can be used also to the following effect (let us recall that a set is punctiform if it contains no non-trivial continuum, cf. \cite{E2}).

\medskip

\noindent {\bf Proposition 5.4.1.}
{\it Let $\alpha \geq 0$ be a countable ordinal and let  $K$ be a compact space with ind$K = \alpha +2$. Then there exists a punctiform $G_{\delta}$-set $H$ in $K$ such that each perfect image of its complement $C = K \setminus H$   contains a compact set whose transfinite dimension ind is either $\geq \alpha$ or is not defined.}

\medskip

In particular, for each natural $n > 1$ there is a finite-dimensional $\sigma$-compact space  whose all perfect images are at least $n$-dimensional.
In fact, the case of natural $n \geq 1$ can be treated separately, without appealing to the Hurewicz theorem, leading to a conclusion that each $(n+1)$-dimensional compact space contains an $n$-dimensional $\sigma$-compact space, all whose perfect images are at least $n$-dimensional.

Let us sketch this reasoning.

Let $K$ be an $(n+1)$dimensional compact space,  $n \geq 1$ being a natural number.
Let $A_{i}$, $B_{i}$, $i=1, \ldots , n+1$, be an essential family of pairs of disjoint closed sets in $K$, cf. \cite{vM2}, sec. 3.1, and let $A_{i}^{\star} \supset A_{i}$, $B_{i}^{\star} \supset B_{i}$ be disjoint closed sets in $K$ containing in its interiors $A_{i}$ and $B_{i}$, respectively.

Using a theorem of Hurewicz, \cite{Ku}, \S 45, IX, one can find continuous maps $p_{i} : K \to [0,1]$ such that $A_{i}^{\star} \subset p_{i}^{-1}(0)$, $B_{i}^{\star} \subset p_{i}^{-1}(1)$ and for some Cantor sets $T_{i} \subset (0,1)$, the fibers $p_{i}^{-1}(t)$ with $t \in T_{i}$ are at most $n$-dimensional, $i = 1, \ldots ,n$.

One can pick $G_{\delta}$-sets $E_{i} \subset p_{i}^{-1}(T_{i})$ such that $E_{i}$ hits every continuum joining $A_{i}^{\star}$ and $B_{i}^{\star}$ and $p_{i}$ is injective on $E_{i}$, cf. \cite{vM2}, proofs of 3.9.3 and 3.11.8. Then, cf. \cite{PP2}, Lemma 5.1, $H_{0} =  E_{1} \cup \ldots \cup E_{n}$ contains no non-trivial  continuum.

The closed set $S = \bigcup _{i=1}^{n} p_{i}^{-1}(T_{i})$ contains $H_{0}$ and is at most $n$-dimensional. Let $H_{1}\subset K \setminus S$ be a zero-dimensional $G_{\delta}$-set such that $(K \setminus S) \setminus H_{1}$ is at most $n$-dimensional, and let $H=H_{0} \cup H_{1}$.

Then $H$ contains no non-trivial continuum and hits every continuum joining some pair $A_{i}^{\star}$, $B_{i}^{\star}$, for $i=1, \ldots , n$.

We claim that the $n$-dimensional $\sigma$-compact set $C = K \setminus H$ has required properties.

Indeed, let $f : C \to L$ be a perfect surjection, and let $u : C \to M$, $v : M \to L$ be the monotone-light decomposition of $f$, i.e., $f = v \circ u$, the fibers of $u$ are connected, the fibers of $v$ are zero-dimensional, and both $u, v$ are perfect, cf. \cite{E1}, 6.2.22. Then $u(A_{i}^{\star} \cap C) \cap u(B_{i}^{\star} \cap C) = \emptyset$, for $i=1, \ldots , n$. Let $S_{i}$ be a partition in $M$ between $u(A_{i}^{\star} \cap C)$ and $u(B_{i}^{\star} \cap C)$ and let $C_{i} = u^{-1}(S_{i})$. Since $C_{i}$ is a partition in $C$ between $A_{i}^{\star} \cap C$ and $B_{i}^{\star} \cap C$
and $A_{i}$ and $B_{i}$ are in the interior of $A_{i}^{\star}$ and $B_{i}^{\star}$, respectively, there exists a partition $K_{i}$ in $K$ between $A_{i}$ and $B_{i}$ such that $K_{i} \cap C \subset C_{i}$, cf. \cite{E2}, Lemma 1.2.9. The intersection $K_{1} \cap \ldots \cap K_{n}$ contains a continuum joining $A_{n+1}$ and $B_{n+1}$, and since $H$ contains no non-trivial continuum,
$ K_{1} \cap \ldots \cap K_{n} \cap C \neq \emptyset$, hence $C_{1} \cap \ldots \cap C_{n} \neq \emptyset$. It follows that $S_{1} \cap \ldots \cap S_{n}\neq \emptyset$ and this shows that $M$ is at least $n$-dimensional.
Since light perfect mappings do not lower the dimension, also $L$ is at least $n$-dimensional.

\subsection{Perfect maps on certain sums of subspaces.}

Let $G$ be an open set in a space $X$ such that $F = X \setminus G$ is locally compact. If $G$ has a perfect map onto a space with ind equal to $\alpha$, then $X$ has a perfect map onto a space with ind $\leq \alpha + 2$.

To see this, let us consider a compact extension $Z$ of $X$ and let $d$ be a metric on $Z$ bounded by 1. Let $\overline{F}$ be the closure of $F$ in $Z$ and $K = \overline{F} \setminus F$ (notice that $K$ is compact, $F$ being locally compact).

Let $f : G \to Y$ be a perfect map onto a space with ind$Y = \alpha$. Let $C(Y)$ be the metric cone over $Y$, i.e., $C(Y) = \{ v \} \cup \,(Y \times (0,1])$, where basic neighbourhoods of the vertex $v$ are the sets $\{ v \} \cup\, ( Y \times (0,\epsilon ))$ and let $E(Y)=(C(Y) \times [0,1]) \setminus \{ (v,0) \}$. We have ind$E(X) \leq \alpha +2$, cf. \cite{E2}, 7.2.F.

The function $g : X \to E(Y) $ defined by $g(x) = \left( f(x), dist (x,\overline{F}), dist (x,K)\right)$ if $x \in G$ and $g(x) = (v, dist (x,K))$ if $x \in F$ is a perfect map (following \cite{E1}, \cite{E2} we assume that the image of a perfect map  is closed in its range).

To check that $g$ is perfect, let us consider $g(x_{n}) \to u$, $u \in E(Y)$. If $u \in (Y \times (0,1]) \times [0,1]$, $dist(x_{n}, \overline{F}) \geq \delta >0$ for some $\delta$ and  almost all $n$, and therefore $x_{n} \in G$ for almost all $n$, and  $f(x_{n})$ converges to the first coordinate of $u$, belonging to $Y$. By perfectness of $f$, the sequence $(x_{n})_{n}$ has a subsequence convergent in $G$, hence in $X$. If $u = (v,c)$, $c>0$, let us pick a subsequence $(x_{n_{k}})_{k}$ of $(x_{n})_{n}$ converging to a point $z$ in $Z$. Since $dist(x_{n}, \overline{F}) \to 0$, $z \in \overline{F}$ and since $dist(x_{n}, K) \to c >0$, $dist(z,K) >0$. In effect, $z \in \overline{F} \setminus K = F$, i.e. the sequence $(x_{n_{k}})_{k}$ is convergent in $X$.

Let us notice the following consequence of this approach.

Let $\D _{n}$ be the class of spaces which can be exhausted in $n$ steps by a subsequent removing of maximal relatively open locally compact sets. Then each element of $\D _{n}$ has a perfect map onto at most $2n-1$-dimensional space. Let $d(n) \leq 2n-1$ be the minimal number with this property. We do not know if the numbers $d(n)$, $n=1,2,\ldots$, are bounded.

To conclude, it seems that the remark of Isbell \cite{Is}, page 119, that the relations between perfect mappings and dimensional properties call for more studies is still valid.

\end{document}